\documentclass[12pt,twoside]{amsart}
\usepackage{amssymb,amsmath,amsthm}
\usepackage{amsfonts}
\usepackage{mathrsfs}
\usepackage{verbatim}
\usepackage{graphicx}
\usepackage{epsfig, enumerate}
\usepackage{soul}
\usepackage{ dsfont }
\usepackage[all]{xy}
\usepackage{color}
\DeclareMathAlphabet{\mathpzc}{OT1}{pzc}{m}{it}

\def\A{\mathcal A}
\def\B{\mathcal B}
\def\H{\mathcal H}
\def\K{\mathcal K}
\def\L{\mathcal L}
\def\P{\mathcal P}
\def\S{\mathcal S}
\def\QN{\mathcal{QN}}
\def\eell{\overline{\ell}}
\def\d{\displaystyle}
\def\ep{\varepsilon}
\def\m{\mathpzc{m}}
\def\pco{\text{{\it p}{\rm -co}}}
\def\qco{\text{{\it q}{\rm -co}}}
\def\coe{{\rm coe-}}
\def\Lip{{\rm Lip}}
\def\Ae{\mbox{\AE}}

\newtheorem{theorem}{Theorem}[section]
\newtheorem{lemma}[theorem]{Lemma}
\newtheorem{definition}[theorem]{Definition}
\newtheorem{corollary}[theorem]{Corollary}
\newtheorem{proposition}[theorem]{Proposition}

\begin{document}

\title{Spaceability of sets of \lowercase{$p$}-compact maps}

\author{Thiago R. Alves}
\thanks{Thiago R. Alves was financed in part by the Coordenação de Aperfeiçoamento de Pessoal de Nível Superior - Brasil (CAPES) - Finance Code 001 and FAPEAM}
\address{Thiago R. Alves\\
	Departamento de matem\'{a}tica \\ Instituto de Ci\^{e}ncias Exatas \\ Universidade Federal do Amazonas \\ 69.077-000 -- Manaus -- Brasil.}
\email{alves@ufam.edu.br}

\author{Pablo Turco}
\thanks{Pablo Turco was supported in part by ANPCyT PICT-2019-00080.}
\address{Pablo Turco\\
	IMAS--UBA--CONICET \\ Pab. I, Facultad de Ciencias Exactas y Naturales, Universidad de Buenos Aires \\ 1428 Buenos Aires. Argentina}
\email{paturco@dm.uba.ar}


\begin{abstract} 
We provide quite sufficient conditions on the Banach spaces $E$ and $F$ in order to obtain the spaceability of the set of all linear operators from $E$ into $F$ which are $q$-compact but not $p$-compact. Also, under similar conditions over $E$, we prove that this set contains (up to the null operator) a copy of $\ell_s$ whenever $F = \ell_s$. Finally, we give some applications of our previous results to show the spaceability of some sets formed by non-linear mappings (polynomial and Lipschitz) which are $q$-compact but not $p$-compact. The spaceability in the space of holomorphic mappings determined by $p$-compact sets is also considered.
\end{abstract}
\maketitle

\section{Introduction}

A well-known result due to Grothendieck states that a subset in a normed space is a relatively compact set if and only if it is contained in the convex hull of some norm null sequence of vectors in the space. Inspired by this characterization, Sinha and Karn \cite{SiKa} have introduced the notion of relatively $p$-compact sets ($1\leq p<\infty$) of a normed space in a similar way. Indeed, they have defined a relatively $p$-compact set via a $p$-convex hull of some norm $p$-summable sequence (see definition below). 
 
Once the notion of $p$-compact sets is established, it becomes natural to introduce the notion of $p$-compact operators in the same way as we define the compact operators. Namely, a linear operator between two Banach spaces is $p$-compact if it maps the closed unit ball into a relatively $p$-compact set. The set of all $p$-compact operators from a Banach space $E$ into a Banach space $F$ is denoted by $\K_p(E,F)$. Several researchers have been developing the theory concerning this class of operators (see, e.g., \cite{DelPinSer, GaLaTur, Pietsch_p_compact, DelPin, SiKa} and the references therein). In the same line, many other classes of mappings were introduced by taking the natural definition of compact mappings and replacing {\it the compact set part} of the definition by $p$-compact sets. For instance, we have the $p$-compact polynomials (see, e.g., \cite{AMR, AR, LaTur1} and the references therein), $p$-compact Lipschitz mappings (see \cite{AchDahTur}) and $p$-compact holomorphic mappings (see, e.g., \cite{ACGM_Trans, AMR, LaTur1}).

Relaying linear $p$-compact operators, in different articles (see, e.g., \cite{ACGM_Trans, GaLaTur, Pietsch_p_compact, DelPin}) are established some assumptions over Banach spaces $E$ and $F$ for which, for some $1\leq p<q<\infty$, $p$-compact and $q$-compact linear operators from $E$ to $F$ coincides or not. The first  main purpose of this article is to investigate the conditions for which there exist $q$-compact mappings (linear, polynomial, Lipschitz, holomorphic) which are not $p$-compact. Moreover, another major goal is studying if there exist (up to the null mapping) linear structures formed by $q$-compact mappings which are not $p$-compact.

This leds us to the concept of lineability and spaceability. More precisely, a subset $S$ of a topological vector space $X$ is called {\it lineable} if there is an infinite dimensional subspace $W\subset X$ such that $W\subset S \cup \{0\}$. If the subspace $W$ is closed, then $S$ is called {\it spaceable}. A seminal example of this kind was provided in \cite{Gu91}. Namely, the family of continuous functions on the closed unit interval $[0,1]$ that are differentiable at no point contains, except for the null function, an infinite dimensional vector space. Besides these kind of results have been established at least since 60's, the notions of lineability and spaceability were formalized just in \cite{AGS, GurQua}. Henceforth, many authors have been working on the subject (see, e.g., \cite{AronBerPelSeo, BerPelSeo} and references therein).

It is worth mentioning that Hern\'andez, Ruiz and S\'anchez \cite[Theorem 3.5]{HeRuSan} have studied the spaceability of sets formed by differences of general Banach operator ideals. In particular, we can deduce from \cite[Theorem 3.5]{HeRuSan} that, under some hypotheses on the Banach spaces $E$ and $F$, the set $\K_q(E,F) \setminus \K_p(E,F)$ is spaceable whenever $\K_p(E,F) \not= \K_q(E,F)$ with $1 \leq p < q$. In this direction, we are able to show that the previous set, for $F = \ell_s$, actually contains copies of $\ell_s$. Also, we study the set $\K_q(E,F) \setminus \K_p(E,F)$ assuming local geometric structures of $E$ and $F$ which are different from those considered in \cite{HeRuSan}.

This article is organized as follows. Section \ref{Sec: First cases}, after fixing some notations and definitions, we study the spaceability of the set formed by $q$-compact linear operators which are not $p$-compact where the domain and codomain are $\ell_r$ and $c_0$. As a first result we show that if the set $\K_q(\ell_r,\ell_s)\setminus \K_p(\ell_r,\ell_s)$ is non-empty, then we may find (up to $0$) a subspace isomorphic to $\ell_s$ inside $\K_q(\ell_r,\ell_s)\setminus \K_p(\ell_r,\ell_s)$ (Theorem~\ref{Thm:subspace lr}). Then, since there exists a relation between $p$-compact operators and $p$-summing operators (Proposition~\ref{Prop: relations}) using some well-known results about when $q$-summing and $p$-summing operators coincide, we exemplify our results with a list in which cases give $\K_q(\ell_r,\ell_s)\setminus \K_p(\ell_r,\ell_s)$ is empty or not.

In Section \ref{Sec: General cases} we extend some results of Section~\ref{Sec: First cases}. In particular, we investigate the existence of linear structures into the set $\K_q(E,F)\setminus \K_p(E,F)$ for more general Banach spaces $E$ and $F$. Nevertheless, we assume that the geometries of $E$ and $F$ are, somehow, similar to $\ell_r$. It is worth noting that the main result of this section (Theorem \ref{Thm: Main_Linear}) depends if $\K_q(\ell_r,\ell_s)\setminus \K_p(\ell_r,\ell_s)$ is empty or not, which we have analyzed in Section~\ref{Sec: First cases}.

In Section~\ref{Sec: non linear case} we investigate the existence of linear structures inside sets of non-linear mappings which are determined by $p$-compact sets. Specifically, we use some results from previous sections in order to study sets of polynomial and Lipschitz mappings which are, each of them, $q$-compact but not $p$-compact. We compare the results obtained and show that, in some sense, there is {\it much more variety of} $q$-compact Lipschitz or polynomials than $q$-compact linear operators. For instance, meanwhile for every $2\leq p<q$, $q$-compact linear operators from $\ell_2$ to $\ell_2$ are $p$-compact, for any $2\leq p<q$ there is a subspace isomorphic to $\ell_2$ which consist (up to $0$) in $q$-compact polynomials which are not $p$-compact (Proposition~\ref{Prop: Comparation}). This section finish with a brief discussion about the spaceability in the space of $p$-compact holomorphic mappings. In \cite{LaTur1} it is shown that, there exists a non $p$-compact holomorphic mapping from $\ell_1$ to $\ell_p$ such that every $n$-homogeneous polynomial of its Taylor series expansion at every point is $p$-compact, giving an answer to a problem posted in \cite{AMR}. Now we show that there exists a closed subspace of this last type of holomorphic function which are not $p$-compact, when we consider the Nachbin topology (Proposition~\ref{Prop: Space_holo}).

The notation we use is standard. The letters $E$ and $F$ will always denote Banach spaces over the same field $\mathbb K=\mathbb R$ or $\mathbb C$, except we deal with holomorphic functions, in which $\mathbb K=\mathbb C$. The open unit ball of $E$ is $B_E$ and the dual space of $E$ is denoted by $E'$. For $1\leq p<\infty$, the space of all $p$-summable sequences of $E$ is denoted by $\ell_p(E)$ and it is a Banach space endowed with the norm $\|(x_n)_n\|_{\ell_p(E)}=(\sum_{n=1}^{\infty}\|x_n\|^p)^{1/p}$. A set $K\subset E$ is called {\it relatively $p$-compact} (in the sense of Sinha and Karn \cite{SiKa}) if there exists a $p$-summable sequence $(x_n)_n$ in $E$ such that 
$$
K\subset \pco\{(x_n)_n\} := \{\sum_{n=1}^{\infty} \alpha_n x_n \colon (\alpha_n)_n \in B_{\ell_{p'}}\},
$$
where $p'$ is the conjugated index of $p$. The {\it $p$-compact measure} of $K$ was introduced in \cite{LaTur1}  as
$$
\m_p(K;E)=\inf\{\|(x_n)_n\|_{\ell_p(E)}\}
$$
where the infimum is taken over all the $p$-summable sequences $(x_n)_n\subset E$ such that $K\subset \pco\{(x_n)_n\}$. The set of $p$-compact operators forms a Banach operator ideal endowed with the norm given by $\|T\|_{\K_p}=\m_p(T(B_E);F)$ for every $T\in \K_p(E,F)$.
It is clear that if $1\leq p<q$, then every $p$-compact set is $q$-compact. Then, we always have the inclusion $\K_p(E,F)\subset \K_q(E,F)$. For a general background of Banach operator ideals we refer the reader to the book of Pietsch \cite{Pietsch_Book}.

Since in most of our results we will deal with $c_0$ instead of $\ell_{\infty}$, we will use the following notation. For $1\leq r \leq \infty$,  
$$
\eell_r=\left\{\begin{array}{cc} \ell_r & \mbox{if\ } r<\infty;\\
	c_0 & \mbox{if\ } r=\infty.
\end{array}\right.
$$ 
Also, $\ell_r^{n}$ stands for the $n$-dimensional vector space endowed with the $\|\cdot\|_{\ell_r}$-norm. All other terminology, definitions and references are in their respective sections.

\section{Case $\K_q(\eell_r,\eell_s)\setminus \K_p(\eell_r,\eell_s)$}\label{Sec: First cases}

First we will dedicate to study the case whether $\K_q(\eell_r,\eell_s)\setminus \K_p(\eell_r,\eell_s)$ is empty or not. The results obtained will be generalized in the next section. From \cite[Theorem~3.5]{HeRuSan} we know that, for any Banach space $E$, the set $\K_q(E,\eell_s)\setminus \K_p(E,\eell_s)$ is empty or spaceable. We will refine this result by showing that in the case of $\K_q(E,\eell_s)\setminus \K_p(E,\eell_s)$ is non-empty, it contains (up to the null operator) a copy of $\eell_s$. It is clear that if a set contains (up to the null vector) a copy of $\eell_s$ then it is spaceable; but the reverse is, in general, not true. First, we introduce a definition.

\begin{definition} \rm \label{S-spaceability} Let $E$ be a Banach space, $M$ be a subset of $E$, and $\S$ be a Banach space. We say that $M$ is {\it $\S$-spaceable} if there exists a copy of $\S$ in $M\cup \{0\}$. Moreover, we say that $M$ is isometrically $\S$-spaceable if there is an isometry from $\S$ into $M\cup\{0\}$.
\end{definition}

The proof of the above mentioned result also allows us to estimate the Banach-Mazur distance between $\ell_s$ and its copy in $\K_q(E,\eell_s)\setminus\K_p(E,\eell_s)$. Recall that the {\it Banach-Mazur distance} between two Banach spaces $E$ and $F$ is given by
$$d(E,F) := \inf\{\|T\| \|T\|^{-1}: T \in \mathcal{L}(E,F) \mbox{ is an isomorphism}\}.$$
Then, $d(E,F)$ is the smallest possible $r \ge 1$ for which there exists an isomorphism $T : E \rightarrow F$ such that
$$B_F \subset T(B_E) \subset r B_F.$$

Let us introduce two operators that we will use in the following two results. Fix $r \in [1,\infty]$ and let $A = \{m_1 < m_2 < m_3 < \cdots \}$ be a subset of $\mathbb N$. We define (by continuous linear extention) two continuous linear operators as follows: 
	\begin{align*}
	\iota_A &: \eell_r \to \eell_r; \ \ \iota_{A}(e_k) :=e_{m_k}; \\ P_A &: \eell_r \to \eell_r; \ \ P_{A}(e_k) :=\left\{\begin{array}{rl}
	e_j &\mbox{if} \ k=m_j\\
	0 &\mbox{otherwise}
\end{array}\right.;
\end{align*}
for every $k \in \mathbb N$.

The proof of the following lemma is straightforward and we will omit it.

\begin{lemma}\label{Lemma: Split ellr}
Let $(\sigma_n)_n$ be a disjoint partition of $\mathbb N$ of infinite subsets, i.e., each $\sigma_n$ is infinite, $\bigcup_n \sigma_n = \mathbb N$ and $\sigma_n \cap \sigma_m = \emptyset$ if $n \not= m$. Then, for $1 \leq r \leq \infty$ and $\lambda = (\lambda_n)_n \in \eell_r$, the linear operator
\begin{align*} 
U_\lambda\colon \eell_r \rightarrow \eell_r; \ \ U_{\lambda}(x) := \sum_{n=1}^{\infty} \lambda_n \iota_{\sigma_n}(x)
\end{align*}
satisfies the following properties:
\begin{enumerate}[\upshape a)]
	\item $U_{\lambda}$ is a $\|\lambda\|_{\eell_r}$-isometry. This is $\|U_{\lambda}(\alpha)\|=\|\lambda\|_{\eell_r}\|\alpha\|_{\eell_r}$ for every $\alpha \in \eell_r$.
	\item $P_{\sigma_n} \circ U_{\lambda}=\lambda_n Id_{\eell_r}$.
\end{enumerate}
\end{lemma}
 
\begin{theorem}\label{Thm:subspace lr} Let $E$ be a Banach space, let $1\leq p<q\leq \infty$, and let $1\leq r\leq \infty$. Then the set $\mathscr{A} := \K_q(E,\eell_r) \setminus \K_p(E,\eell_r)$ is either empty or $\eell_r$-spaceable. If $\mathscr{A} \not= \emptyset$, then there is a subspace $M \subset \mathscr{A} \cup \{0\}$ isomorphic to $\eell_r$ satisfying
\begin{align*}
	d(M,\eell_r)\leq \inf_{S \in \mathscr{A}} \dfrac{\|S\|_{\K_q}}{\|S\|}.
\end{align*}
Moreover, for $r\geq q$, the set $\mathscr{A}$
is isometrically $\ell_r$-spaceable.
\end{theorem}
\begin{proof}
Suppose there is $S \in \mathscr{A}$. By using the same notations of Lemma \ref{Lemma: Split ellr}, we define an operator $\Psi\colon \eell_r\rightarrow \K_q(E,\eell_r)$ such that
$$\Psi(\lambda) = U_{\lambda}\circ S.$$
It is clear that $\Psi$ is a well-defined linear operator and, since
\begin{align} \label{Aux-Norm-Psi-1}
\|\Psi(\lambda)\|_{\K_q}\leq \|U_{\lambda}\| \|S\|_{\K_q}=\|\lambda\|_{\eell_r}\|S\|_{\K_q},
\end{align}
$\Psi$ is continuous with $\|\Psi\| \leq \|S\|_{\K_q}$. On the other hand, since $U_\lambda$ is a $\|\lambda\|_{\eell_r}$-isometry (by Lemma \ref{Lemma: Split ellr}), we have
\begin{align} \label{Aux-Norm-Psi-2}
\|\Psi(\lambda)\|_{\K_q} \geq \|\Psi(\lambda)\| = \sup_{x \in B_E} \|U_\lambda(S(x))\|_{\eell_r} = \|\lambda\|_{\eell_r}\|S\|.
\end{align}
It follows from (\ref{Aux-Norm-Psi-1}) and (\ref{Aux-Norm-Psi-2}) that $\Psi$ is an isomorphism from $\eell_r$ onto $\Psi(\eell_r)$ and the norm of its inverse $\Psi^{-1}$  satisfies
$$\|\Psi^{-1}\| \leq \dfrac{1}{\|S\|}.$$
Putting this all together we get
\begin{equation}\label{eq:psinorm}
\|S\|\leq \|\Psi\|\leq \|S\|_{\K_q} \ \ \mbox{and} \ \ d(\Psi(\eell_r), \eell_r) \leq \dfrac{\|S\|_{\K_q}}{\|S\|}.
\end{equation}
To finish the first part of the proof it remains to show that $\Psi(\lambda)$ is not a $p$-compact operator for any $0\neq\lambda \in \eell_r$. To prove this, we take $\lambda = (\lambda_n) \in \eell_p$ with $\lambda_{n_0} \not= 0$ for some $n_0 \in \mathbb N$. It follows from Lemma \ref{Lemma: Split ellr}~b) that $S=\frac{1}{\lambda_n}P_{\sigma_n}\Psi(\lambda)$, which implies (by the ideal property for $p$-compact operators) that $\Psi(\lambda)$ is not a $p$-compact operator.

For the moreover part, we will show the case when $r<\infty$. The case $r=\infty$ follows in the same way. Fix $\lambda\in \ell_r$. As $\Psi(\lambda)$ is $q$-compact, there exists a sequence $(z_j)_j \in \ell_q(\ell_r)$ such that $\Psi(\lambda)(B_E)\subset \qco\{(z_j)_j\}$. Hence, for each $n \in \mathbb N$, we have
$$
\lambda_n S(B_E) = P_{\sigma_n}\Psi(\lambda)(B_E) \subset \qco\{(P_{\sigma_n}(z_j))_j\},
$$
which implies  
$$\d |\lambda_n|\|S\|_{\K_q} \leq \|(P_{\sigma_n}(z_j))_j\|_{\ell_q(\ell_r)}=\left(\sum_{j=1}^{\infty}\left(\sum_{k\in \sigma_n}|e'_k(z_j)|^r\right)^{q/r}\right)^{1/q},$$
where $e'_k$ is the coordinate functional. Then, since we are assuming that $r\geq q$, we may apply the Minkowski inequality to obtain
$$
\begin{array}{rl}
\d \|S\|^r_{\K_q}\sum_{n=1}^{\infty} |\lambda_n|^r\leq& \d \sum_{n=1}^{\infty} \left(\sum_{j=1}^{\infty}\left(\sum_{k\in \sigma_n}|e'_k(z_j)|^r\right)^{q/r}\right)^{r/q}\\
\d =&\d \left(\sum_{n=1}^{\infty} \left(\sum_{j=1}^{\infty}\left(\sum_{k\in \sigma_n}|e'_k(z_j)|^r\right)^{q/r}\right)^{r/q}\right)^{q/r \ r/q} \\
\d \leq&\d \left(\sum_{j=1}^{\infty}\left(\sum_{n=1}^{\infty} \sum_{k\in \sigma_n}|e'_k(z_j)|^r\right)^{q/r}\right)^{r/q}\\
\d =&\d \left(\sum_{j=1}^{\infty} \left(\sum_{k=1}^{\infty} |e'_k(z_j)|^r\right)^{q/r}\right)^{r/q}\\
\d =&\d \|(z_j)_j\|_{\ell_q(\ell_r)}^r.
\end{array}
$$
Note that we have just proven that $\|S\|_{\K_q} \|\lambda\|_{\eell_r} \leq \|(z_j)_j\|_{\ell_q(\ell_r)}$ whenever $\Psi(\lambda) \subset \qco\{(z_j)_j\}$.  From this, we get
$$\|S\|_{\K_q} \|\lambda\|_{\ell_r} \leq \inf\{\|(z_j)_j\|_{\ell_q(\ell_r)}\colon \Psi(\lambda) \subset \qco\{(z_j)_j\}\}=\|\Psi(\lambda)\|_{\K_q}.$$ 
This inequality together with inequality \eqref{Aux-Norm-Psi-1} give us
$$
\|S\|_{\K_q} \|\lambda\|_{\ell_r} = \|\Psi(\lambda)\|_{\K_q}.
$$
This is, $\Psi$ is a $\|S\|_{\K_q}$-isometry. So the proof is completed by taking $S \in \mathscr{A}$ such that $\|S\|_{\K_q} = 1$.
\end{proof}

From \cite[Theorem~3.13]{GaLaTur} follows that $\K_q(E,F) = \K_p(E,F)$ if the spaces of $q$-summing and  $p$-summing operators defined from $F'$ into $E'$ coincide. As we will see bellow, for some cases we actually have an equivalence. To prove this we need some terminology and a lemma. Recall that for a Banach operator ideal $\A$, the {\it dual ideal} of $\A$, denoted by $\A^d$, is defined for Banach spaces $E$ and $F$ as
$$
\A^{d}(E,F) := \{T\in \L(E,F)\colon T'\in \A(F',E')\}
$$
with the norm $\|T\|_{\A^d} :=\|T'\|_\A$, where $T'$ represents the adjoint of $T$.
\begin{lemma}\label{Lemma: Dual ideals}
Let $\A$ and $\B$ be Banach operator ideals, and let $E$ and $F$ be Banach spaces. Suppose that $\A=\A^{dd}$, then, if $\B^{d}(E,F')\setminus \A^{d}(E,F')$ is non-empty, $ \B(F,E')\setminus \A(F,E')$ is non-empty. 
\end{lemma}
\begin{proof}
Take $T\in  \B^{d}(E,F')\setminus \A^{d}(E,F')$ and consider $S=T'\circ J_F$, where $J_F\colon F\rightarrow F''$ is the canonical inclusion. Then $S\in \B(F,E')$ and, if $S$ belongs to $\A(F,E')$, since $\A=\A^{dd}$ we have that $T=S'\circ J_E \in \A^d(E,F')$, which is a contradiction.
\end{proof}

For $1\leq p <\infty$, we let $\Pi_p$ denote the Banach operator ideal of  all $p$-summing linear operators. Meanwhile the Banach operator ideal of quasi $p$-nuclear operators is denoted by $\QN_p$. From \cite[Theorem 2.8]{GaLaTur} we have that $\K_p= \QN_p^d$ and that $\K_p^d = \QN_p$. Also, it follows from \cite[Theorem 10]{Pietsch_p_compact} that the maximal hull of the ideal of quasi $p$-nuclear operators is the ideal of $p$-summing operators: $\QN_p^{max}=\Pi_p$. From this and \cite[4.9.12]{Pietsch_Book} we get that the minimal kernel of both ideals coincide: $\QN_p^{min}=\Pi_p^{min}$.

\begin{proposition}\label{Prop: relations}
Let $E$ be a Banach space, $1\leq p<q<\infty$ and $1\leq r< \infty$. Then
\begin{enumerate}[\upshape a)]
\item $\K_q(E,\ell_r)\setminus \K_p(E,\ell_r)$ is empty if and only if $\Pi_q(\ell_{r'},E')\setminus \Pi_p(\ell_{r'},E')$  is empty ($\ell_{r'}=c_0$ if $r=1$).
\item $\K_q(E,c_0)\setminus\K_p(E,c_0)$ is empty if and only if $\Pi_q(\ell_1,E')\setminus \Pi_p(\ell_1,E')$ is empty.
\end{enumerate}
\end{proposition}
\begin{proof}
For item a), we first note that by combining \cite[Corollary~1]{Per} and \cite[Theorem 5]{Per} we have that $\Pi_q(\ell_{r'},E')$ coincides with $\QN_q(\ell_{r'}, E')$ (see the last paragraph of \cite{Per}).  This gives the equality $\K_q(E,\ell_r)=\Pi^d_q(E,\ell_r)$. Also, as $\K_p^{dd}=\K_p$ \cite[Proposition~8]{Pietsch_p_compact} and $\QN_p^{dd} = \QN_p$, the proof follows by an application of Lemma~\ref{Lemma: Dual ideals}. 

For item b), for the {\it if} part, take $T \in \K_q(E,c_0)$ and since $c_0$ has the approximation property, there exists a sequence of finite rank operators $(T_n)_n$ which converges in the $\|\cdot\|_{\K_q}$ norm and, by \cite[Proposition~2.9]{GaLaTur}, coincides with the $\|\cdot\|_{\Pi^d_q}$ norm. Thus, $\|T_n\|_{\K_q}=\|T_n'\|_{\Pi_q}$ and since  $\Pi_q(\ell_1,E')=\Pi_p(\ell_1,E')$, by the closed graph theorem, we get a constant $C>0$ such that $\|T_n'\|_{\Pi_q}\leq \|T_n'\|_{\Pi_p}\leq C \|T'_n\|_{\Pi_q}$. This implies that the sequence of finite rank operators $(T_n)_n$ also converges to $T$ in the $\|\cdot\|_{\K_p}$ norm and $T$ is a $p$-compact operator. 

For the {\it only if part}, since $\K_q(E,c_0)=\K_p(E,c_0)$ there exists $C>0$ such that $\|\cdot\|_{\K_q}\leq \|\cdot\|_{\K_p}\leq C\|\cdot\|_{\K_q}$ and to see that $\Pi_q(\ell_1,E')$ coincides with $\Pi_p(\ell_1,E')$ is equivalent to see that $\QN_q(\ell_1,E')$ coincides with $\QN_p(\ell_1,E')$. This is because since $(\ell_1)'$ has the approximation property, $\QN_p(\ell_1,E')=\QN_p^{min}(\ell_1,F')$ and by \cite[Proposition~22.8]{DF} $\QN_p^{min}=\QN_q^{min}$ if and only if $\Pi_p=\QN_p^{max}=\QN_q^{max}=\Pi_q$. Now, take $T\in \QN_q(\ell_1,E')$ and applying \cite[Corollary~2.4]{FouSwa} we may find a compact operator $R$ from $\ell_1$ to $\ell_1$ and a quasi $q$-nuclear operator $S$ from $\ell_1$ to $E'$ such that $T=S\circ R$. Consider $P_n\colon \ell_1\rightarrow \ell_1$ the projection of $\ell_1$ onto the first $n$-coordinate and define, for each $n\in \mathbb N$ the finite rank operator $T_n=S \circ \pi_n\circ R$. Routing arguments show that $T_n$ converges in the $\QN_q$-norm to $T$, so the proof will follows by showing that $\|T_n\|_{\QN_q}\leq \|T_n\|_{\QN_p}\leq C \|T_n\|_{\QN_q}$ for each $n \in \mathbb N$.

Denoting by $\ell_1^n$ and $\ell_{\infty}^n$ the $n$-dimensional Banach spaces endowed with the norm $\ell_1$ and $\ell_{\infty}$ respectively, it is not difficult to see that
$$
\|S\circ P_n\colon \ell_1\rightarrow E'\|_{\QN_q}=\|S\circ P_n\colon \ell_1^n\rightarrow E'\|_{\QN_q}.
$$
Then, $\|S\circ P_n\colon \ell_1^n\rightarrow E'\|_{\QN_q}=\|(S\circ P_n)'\colon E''\rightarrow \ell^n_{\infty}\|_{\K_q}$ and if $(P_n)'$ is the projection of $\ell_\infty$ to the first $n$-coordinate and $\iota_n\colon \ell_{\infty}^n\rightarrow c_0$ the canonical inclusion, we have that $(S\circ P_n)'=P_n'\circ \iota_n\circ P_n'\circ S'$ and
$$\begin{array}{rl}
\|(S\circ P_n)'\colon E''\rightarrow \ell^n_{\infty}\|_{\K_q}=&\|P_n'\circ \iota_n\circ P_n'\circ S'\colon E''\rightarrow \ell^n_{\infty}\|_{\K_q}\\
\leq& \|\iota_n\circ P_n'\circ S'\colon E''\rightarrow c_0\|_{\K_q}
\\
\leq & C \|\iota_n\circ P_n'\circ S'\colon E''\rightarrow c_0\|_{\K_p}\\
\leq & C \|P_n'\circ S'\colon E''\rightarrow \ell^n_{\infty}\|_{\K_p}.
 \end{array}
$$
This give us that $\|T_n\|_{\QN_q}\leq \|T_n\|_{\QN_p}\leq C \|T_n\|_{\QN_q}$ as we wanted to see.
\end{proof}

Combining Theorem~\ref{Thm:subspace lr} and Proposition~\ref{Prop: relations}, we may show that in the cases where $\Pi_q(\eell_{r'},E')$ does not coincide with $\Pi_p(\eell_{r'},E')$, there are {\it a lot of} $q$-summing operators which are not $p$-summing.

\begin{theorem}\label{Thm: Space in p-summin}
Let $E$ be a Banach space, let $1\leq p<q\leq \infty$, and let $1\leq r\leq \infty$. 

\begin{enumerate}[\upshape (a)]
\item If $r\neq 1$ and $\Pi_q(\eell_{r'},E')\setminus \Pi_p(\eell_{r'},E')$ is non-empty, then
it is $\eell_{r}$-spaceable. Moreover, there exists $M\subset \Pi_q(\eell_{r'},E')\setminus \Pi_p(\eell_{r'},E')\cup\{0\}$ isomorphic to $\eell_{r}$
with the Banach-Mazur distance between $M$ and $\eell_{r}$ satisfying
$$d(M,\eell_r)\leq \inf\left\{\dfrac{\|S\|_{\Pi_q}}{\|S\|}\colon S\in \Pi_q(\eell_r,E')\setminus \Pi_p(\ell_r,E')\right\}.
$$
Also, if $r\geq q$ then $\Pi_q(\eell_{r'},E')\setminus \Pi_p(\eell_{r'},E')$
is isometrically $\eell_{r}$-spaceable.

\item If $\Pi_q(\ell_1,E')\setminus \Pi_p(\ell_1,E')$ is non-empty, then $\Pi_q(\ell_1,E')\setminus \QN_p(\ell_1,E')$
is $c_0$-spaceable. Moreover, there exists $M\subset \Pi_q(\ell_1,E')\setminus \QN_p(\ell_1,E')\cup\{0\}$ isomorphic to $c_0$
with the Banach-Mazur distance between $M$ and $c_0$ satisfying
$$d(M,c_0)\leq \inf\left\{\dfrac{\|S\|_{\Pi_q}}{\|S\|}\colon S\in \Pi_q(\ell_1,E')\setminus \QN_p(\ell_1,E')\right\}.
$$
\end{enumerate}
\end{theorem}
\begin{proof}
By Theorem~~\ref{Thm:subspace lr} and Proposition~\ref{Prop: relations}, if $\Pi_q(\eell_{r'},E')\setminus \Pi_p(\eell_{r'},E')$ is non-empty, there exists a Banach space ${\widetilde{M}} \subset \K_q(E,\eell_r)\setminus  \K_p(E,\eell_r) \cup \{0\}$ which is isomorphic to $\eell_r$ and $d(\widetilde{M}, \eell_r)$ satisfies the inequality in Theorem \ref{Thm:subspace lr}. Then, the set $$M=\{S\colon \eell_{r'}\rightarrow E' \colon S=T' \quad \mbox{with} \quad T\in \widetilde{M}\}$$ is the Banach space we are looking for. 
\end{proof}

In order to see for which Banach spaces $E$ and $F$ the set $\K_q(E,F)\setminus \K_p(E,F)$, the first approach is to see for which Banach spaces $F$ the $q$-compact sets and the $p$-compact sets coincide or not. In \cite{Pietsch_p_compact} and \cite{DelPin} was established a criterion for Banach spaces in which $q$-compact sets and $p$-compact sets coincide. A similar characterization can be obtained directly from \cite[Proposition~1.8]{LaTur2}.
\begin{proposition}\label{Prop: q=p-compact}
Let $F$ be a Banach space and let $1\leq p<q<\infty$. Then every $q$-compact set is $p$-compact if and only if $\K_q(\ell_1,F)\setminus \K_p(\ell_1,F)$ is empty.
\end{proposition}
\begin{proof}
The only if part is trivial. For the if part, take $K\subset F$ a $q$-compact set. By \cite[Proposition~1.8]{LaTur2}, there exists an operator $T\in \K_q(\ell_1,F)$ such that $K\subset T(B_{\ell_1})$. Since we are assuming that every $q$-compact set from $\ell_1$ to $F$ is $p$-compact, we obtain that $K$ is $p$-compact, and the proof is complete.
\end{proof}

As a consequence we have,

\begin{proposition}\label{Prop: example q=p-comp s>2}
Let $1\leq p<q<\infty$ and $s\geq 2$. Then
\begin{enumerate}
\item If $q\leq 2$, every $q$-compact set in $\eell_s$ is $p$-compact.
\item If $q>2$, there exists a $q$-compact set in $\eell_s$ which is not $p$-compact.
\end{enumerate}
\end{proposition}
\begin{proof}
By the above proposition, we must see if $\K_q(\ell_1,\eell_s)\setminus \K_p(\ell_1,\eell_s)$ is empty or not which, by Proposition~\ref{Prop: relations}, is equivalent to see if $\Pi_q(\ell_{s'},\ell_{\infty})\setminus \Pi_p(\ell_{s'},\ell_{\infty})$ is empty or not. As $1\leq s'\leq 2$, by a Theorem of Kwapi\'en (see \cite[Theorem~2]{Pietsch_p_summ}) if $1\leq s'\leq 2$, then $\Pi_q(\ell_{s'},\ell_{\infty})\setminus \Pi_p(\ell_{s'},\ell_{\infty})$ is empty for $q\leq 2$ and, by \cite[Final Remark, 2]{Pietsch_p_summ}, $\Pi_q(\ell_{s'},\ell_{\infty})\setminus \Pi_p(\ell_{s'},\ell_{\infty})$ is non-empty for $q>2$.
\end{proof}

\begin{proposition}\label{Prop: example q=p-comp s<2}
Let $1\leq p<q<\infty$ and $s\leq 2$. Then
\begin{enumerate}
\item If $q\leq s$, every $q$-compact set in $\ell_s$ is $p$-compact.
\item If $q>s$, there exists a $q$-compact set in $\ell_s$ which is not $p$-compact.
\end{enumerate}
\end{proposition}

\begin{proof} The proof is analogous of the above, using a result of Saphar (see \cite[Theorem~1]{Pietsch_p_summ}) for item $(1)$ and \cite[Final Remark, 3]{Pietsch_p_summ} for item $(2)$.
\end{proof}

We finish the section showing the spaceability (or not) of $\K_q(\eell_r,\eell_s)\setminus \K_p(\eell_r,\eell_s)$ for $r\geq 2$. A first result shows that, despite of in $\eell_s$ for  $s\geq 2$ there exist $q$-compact sets which are not $p$-compact, $q$-compact and $p$-compact operators with range in $\eell_s$ coincide when the domain is $\eell_r$ with $r\geq 2$.

\begin{proposition}\label{Prop: r>2,s>2} Let $1\leq p<q<\infty$ and $r, s\geq 2$. Then $\K_q(\eell_r,\eell_s)\setminus \K_p(\eell_r,\eell_s)$ is empty.
\end{proposition}
\begin{proof}
Fix $r,s\geq 2.$ Using Kwapi\'en's results (see \cite[Theorem~2 and Theorem~3]{Pietsch_p_summ}), for $1\leq p<q<\infty$, we have that $\Pi_q(\eell_{s'},\eell_{r'})\setminus\Pi_p(\eell_{s'},\eell_{r'})$ is an empty set. An application of Proposition~\ref{Prop: relations} completes the proof.
\end{proof}

\begin{proposition}\label{Prop: r>2,s<2} Let $1\leq p<q<\infty$, $r\geq 2$ and $s\leq 2$. Then 
\begin{enumerate}[\upshape a)]
\item If $1<q\leq s$, then  $\K_q(\eell_r,\eell_s)\setminus \K_p(\eell_r,\eell_s)$ is empty.

\item If $s<q\leq 2$, then $\K_q(\eell_r,\eell_s)\setminus \K_p(\eell_r,\eell_s)$ is $\eell_s$-spaceable.

\item If $p<2\leq q$, then $\K_q(\eell_r,\eell_s)\setminus \K_p(\eell_r,\eell_s)$ is $\eell_s$-spaceable.
\item If $p\geq 2$, then $\K_q(\eell_r,\eell_s)\setminus \K_p(\eell_r,\eell_s)$ is empty.
\end{enumerate}
\end{proposition}
\begin{proof}
Item a) follows from Proposition~\ref{Prop: example q=p-comp s<2}. For itens b) and c), for  $s<q\leq 2$ or $p<2\leq q$ it follows from \cite[Final remark]{Pietsch_p_summ} that $\Pi_q(\eell_{s'},\eell_{r'})\setminus \Pi_p(\eell_{s'},\eell_{r'})$ is non-empty. An application of Proposition~\ref{Prop: relations} together with Theorem~\ref{Thm: Space in p-summin} gives the result. Finally, by a result of Kwapi\'en (see \cite[Theorem~3]{Pietsch_p_summ}) for $2\leq p<q$, $r\geq 2$ and $s\leq 2$, we have that $\Pi_q(\eell_{s'},\eell_{r'})\setminus\Pi_p(\eell_{s'},\eell_{r'})$ is an empty set. Another application of Proposition~\ref{Prop: relations} completes the proof.
\end{proof}

\section{Case $\K_q(E,F)\setminus\K_p(E,F)$}\label{Sec: General cases}

In this section we will extend the results obtained about the spaceability of $\K_q(\eell_s,\eell_r)\setminus \K_p(\eell_s,\eell_r)$ to $\K_q(E,F)\setminus\K_p(E,F)$ in the case when the Banach spaces $E$ and $F$ have some part of its geometry, somehow, similar to $\eell_s$ and $\eell_r$. More precisely, we say that a Banach space $E$ contains {\it uniformly copies} of $\ell_r^n$ ($1 \leq r \leq \infty$) if there exists a sequence of subspaces $E_n$ of $E$ such that
\begin{eqnarray} \label{def-embed-unf}
	\sup_{n \in \mathbb N} d(E_n, \ell_r^n) < \infty.
\end{eqnarray}
Equivalently this means that there exists a sequence of injective operators $T_n\colon \ell_r^n\rightarrow E$ such that $\sup_{n\in \mathbb N} \|T_n\| \|T_n^{-1}\| < \infty$. In this case, we will also say that $\ell_r^n$ {\it embed uniformly} in the Banach space $E$.

Also, we say that a Banach space $E$ contains {uniformly complemented copies} of $\ell_r^n$ ($1 \leq r \leq \infty$) if there exist a sequence of subspaces $E_n$ of $E$ and a sequence of projections $P_n \colon E \rightarrow E_n$ satisfying
\begin{eqnarray} \label{def-unf-comp}
\sup_{n \in \mathbb N} d(E_n, \ell_r^n) < \infty \ \ \mbox{ and } \ \ \sup_{n \in \mathbb N} \|P_n\| < \infty.
\end{eqnarray}


The next four lemmas are, in some sense, saying that the $\m_q$-measure of $q$-compact sets in $\eell_r$ is well behaved with respect to the natural projections on finite coordinates. In what follows, $P_n\colon \eell_r\rightarrow \eell_r$ stands for the projection onto the first $n$-coordinates.

\begin{lemma}\label{Lemma: On domain}
Let $1\leq p<q< \infty$, $1\leq r\leq \infty$ and $F$ be a Banach space. Suppose that $\K_q(\eell_r,F)\setminus\K_p(\eell_r,F)$ is non-empty. Then, there exists an operator $R\in \L(\eell_r,F)$ such that, given $\ep>0$ and $M>0$, there exist $n<m \in \mathbb N$  such that 
$$
\|R(P_m-P_n)\|_{\K_q}\leq \ep \quad \mbox{and} \quad \|R(P_m-P_n)\|_{\K_p}>M.
$$
\end{lemma}
\begin{proof}
The first case is for $r=1$. If $\K_q(\ell_1,F)\setminus\K_p(\ell_1,F)$ is not empty, there exists a $q$-compact set in $F$ which is not $p$-compact. Thus, by \cite[Lemma~5]{Tur}, there exists a null sequence $(x_n)_n \subset F$ and an increasing sequence of integer $1 = j_1 < j_2 < j_3 < \cdots$ such that, if $L_m=\{x_{j_m},x_{j_{m}+1},\ldots,x_{j_{m+1}-1}\}$, then
$$
\lim_{m\to \infty} \m_q\left(L_m;F\right)=0 \quad \mbox{and} \quad \lim_{m\to \infty} \m_p\left(L_m;F\right)=\infty.
$$
Define the operator $R\colon \ell_1\to F$ as $R(e_n)=x_n$ for any $n \in \mathbb N$ and extend it by linearity and continuity. For each $m\in \mathbb N$, we have that $R(P_{j_{m+1}-1}-P_{j_m})(B_{\ell_1})=\coe\{L_m\}$, implying that 
$$
\|R(P_{j_{m+1}-1}-P_{j_m})\|_{\K_q}=\m_q\left(L_m;F\right) \quad \mbox{and} \quad \|R(P_{j_{m+1}-1}-P_{j_m})\|_{\K_p}= \m_p\left(L_m;F\right).
$$
For $r\neq 1$, take $T\in \K_q(\eell_r,F)\setminus\K_p(\eell_r,F)$. Since $\K_q=\K_q\circ \K$,  combining \cite[Corollary~2.4 (b) and Theorem 2.5]{FouSwa} we may find a sequence $(\gamma_j)_j\in c_0$ which defines the diagonal operator $D_{(\gamma_j)_j}\colon \eell_r\rightarrow \eell_r$ and an operator $\widetilde T\in \K_q(\eell_r,F)$ such that $T=\widetilde T \circ D_{(\gamma_j)_j}$. We may assume that  $0<\gamma_j<1$ for all $j \in \mathbb N$. 

Define $R\in \L(\eell_r,F)$ as $R=\widetilde T \circ D_{(\sqrt{\gamma_j})_{j}}$. For $n<m$, we have 
$$
\|R(P_m-P_n)\|_{\K_q}\leq \|\widetilde T\|_{\K_q} \|D_{(\sqrt{\gamma_j})_{j}}(P_m-P_n)\|= \|\widetilde T\|_{\K_q} \sup_{n<j\leq m} \{\sqrt{\gamma_j}\}.
$$
It remains to show that $\|R(P_m-P_n)\|_{\K_p}$ is unbounded for every $n,m \in \mathbb N$.  Suppose towards contradiction that there is $M>0$ such that $\|R(P_m-P_n)\|_{\K_p}<M$ for all $n,m \in \mathbb N$. Consider the sequence of operators $(T\circ P_n)_{n \in \mathbb N}$ which converges pointwise to $T$. We claim that the sequence is a $\|\cdot\|_{\K_p}$-Cauchy sequence. Indeed, for $n, m \in \mathbb N$,

$$
\begin{array}{rl}
\|(T\circ P_n)-(T\circ P_m)\|_{\K_p}=&\|\widetilde T \circ D_{(\gamma_j)_j} \circ (P_m-P_n)\|_{\K_p}\\
=&\|R\circ D_{(\sqrt{\gamma_j})_j} \circ (P_m-P_n)\|_{\K_p}\\
=&\|R\circ (P_m-P_n) \circ D_{(\sqrt{\gamma_j})_j} \circ (P_m-P_n)\|_{\K_p}\\
\leq &\|R\circ (P_m-P_n)\|_{\K_p} \|D_{(\sqrt{\gamma_j})_j} \circ (P_m-P_n)\|\\
\leq & \d M \sup_{n<j\leq m} \{\sqrt{\gamma_j}\}.
\end{array}
$$
Thus, $T\in \K_p(\eell_r,F)$, which is a contradiction. So the proof is complete.
\end{proof}

We can reformulate the above lemma in terms of $\ell_r^n$, the Banach space of dimenson $n$ endowed with the $\ell_r$-norm.

\begin{lemma}\label{Lemma: On domain2}
Let $1\leq p<q< \infty$, $1\leq r\leq \infty$ and $F$ be a Banach space. Suppose that $\K_q(\eell_r,F)\setminus\K_p(\eell_r,F)$ is non-empty. Then, for $\ep>0$ and $M>0$, there exist $n \in \mathbb N$ and an operator $R_n\in \L(\ell_r^n,F)$ such than
$$
\|R_n\|_{\K_q}\leq \ep \quad \mbox{and} \quad \|R_n\|_{\K_p}>M.
$$
\end{lemma}

\begin{lemma}\label{Lemma: On Codomain}
Let $1\leq p< \infty$ and $1\leq r\leq \infty$. If $K\subset \eell_r$ is $q$-compact, then 
$\d \lim_{m\to \infty} \m_q(P_m(K);\eell_r)=\m_q(K;\eell_r)$.
\end{lemma}
\begin{proof}
First note that for natural numbers $n<m$, we have $$\m_q(P_n(K);\eell_r) \leq \m_q(P_m(K);\eell_r)\leq \m_q(K;\eell_r).$$ Thus, $\lim\limits_{m\to \infty} \m_q(P_m(K);\eell_r)$ exists. Also, as $K=P_m(K)+(Id-P_m)(K)$ for all $m\in \mathbb N$, we get
$$
\m_q(P_m(K);\eell_r)\leq \m_q(K;\eell_r)\leq \m_q(P_m(K);\eell_r)+\m_q((Id-P_m)(K);\eell_r)
$$
for all $m \in \mathbb N$. If we show that $\lim_{m\to \infty} \m_q((Id-P_m)(K);\eell_r)=0$, we conclude the proof. 

Take a sequence $(x_j)_j \in \ell_q(\eell_r)$ such that $K\subset \qco\{(x_j)_j\}$. Given $\ep>0$, we choose $n,m \in \mathbb N$ such that 
$$\d \sum_{k=n}^{\infty} \|x_k\|_{\eell_r}^p \leq \frac{\ep}{2} \ \ \mbox{ and } \ \ \d \|(Id-P_m)(x_j)\|_{\eell_r}\leq \left(\frac{\ep}{2(n-1)}\right)^{1/p}$$
for each $j=1,\ldots, n-1$. Thus, $(Id-P_m)(K)\subset \qco\{\left((Id-P_m)(x_j)\right)_j\}$, and 
$$
\begin{array}{rl}
\d \m_p((Id-P_m)(K);\eell_r)^p \leq& \d \sum_{j=1}^{\infty} \|(Id-P_m)(x_j)\|_{\eell_r}^{p}\\
=&\d \sum_{j=1}^{n-1} \|(Id-P_m)(x_j)\|_{\eell_r}^{p}+\sum_{j=n}^{\infty} \|(Id-P_m)(x_j)\|_{\eell_r}^{p}\\
\leq& \d \frac{\ep}{2}+\sum_{j=n}^{\infty} \|x_j\|_{\eell_r}^{p}\\
\leq& \d\frac{\ep}{2}+\frac{\ep}{2}=\ep.
\end{array}
$$
As $\ep>0$ was arbitrary, the proof is complete.
\end{proof}

Our last lemma is clear and we omit the proof. For $m\in \mathbb N$, $\widetilde P_m\colon \eell_r\rightarrow \ell^m_r$ denotes the projection of the first $m$-coordinates onto $\ell^m_r$.

\begin{lemma}\label{lemma: mp size finite dimensional}
Let $1\leq q<\infty$, $1\leq r \leq \infty$ and $K\subset \eell_r$ be a bounded set. For any $m \in \mathbb N$, $\m_q(P_m(K);\eell_r)=\m_q(\widetilde P_m(K);\ell^m_r)$.
\end{lemma}

We are almost ready to achieve the main result of this section. For it we will need the following result which provides an important tool  to prove spaceability.

\begin{proposition}[{\cite[Theorem 7.4.2]{AronBerPelSeo}} or {\cite[Theorem~3.3]{KitTim}}] \label{spaceable_codition}
	Let $E_n$ $(n \in \mathbb{N})$ be a Banach space and let $X$ be a Fr\'{e}chet space. Suppose $T_n \in \mathcal{L}(E_n,X)$ for each $n \in \mathbb{N}$ and $Y := span\left\{\bigcup_{n=1}^\infty T_n(E_n)\right\}$. If $Y$ is not closed in $X$ then $X \setminus Y$ is spaceable.
\end{proposition}

\begin{theorem}\label{Thm: Main_Linear} Let $E$ and $F$ be Banach spaces, $1\leq p<q<\infty$ and $1\leq r,s\leq \infty$. Suppose that $\K_q(\eell_r,\eell_s)\setminus \K_p(\eell_r,\eell_s)$ is not empty and that $E$ contains uniformly complemented copies of $\ell_r^n$. If
\begin{enumerate}
\item[a)] $\ell_s^n$ embed uniformly to $F$ and $\K_q(E,F) \setminus \K_p(E,F)$ is not empty,
\end{enumerate}
or  
\begin{enumerate}
\item[b)] $F$ contains uniformly complemented copies of $\ell_s^n$,
\end{enumerate}
then $\K_q(E,F) \setminus \K_p(E,F)$ is spaceable.
\end{theorem}

\begin{proof}
Since the identity $Id\colon \K_p(E,F)\rightarrow \K_q(E,F)$ is continuous  and  $\K_q(E,F) \setminus \K_p(E,F) \not= \emptyset$,  by Proposition~\ref{spaceable_codition} it is sufficient to prove that $\K_p(E,F)$ is not closed in $\K_q(E,F)$. By the closed graph theorem  this follows promptly from
	\begin{eqnarray} \label{eq_mainthm}
	\sup_{T \in \K_p(E,F) \setminus \{0\}} \left({\|T\|_{\K_p}}/{\|T\|_{\K_q}}\right) = \infty.
	\end{eqnarray}
Hence it remains to prove \eqref{eq_mainthm}. 

Take a sequence of positive numbers $(a_j)_j \in \ell_1$ and an increasing sequence of positive numbers $(b_j)_j$ with $\lim\limits_{j\to +\infty} b_j=\infty$. Since $\K_q(\eell_r,\eell_s)\setminus \K_p(\eell_r,\eell_s)$ is non-empty, it follows from Lemma~\ref{Lemma: On domain2} that for each  $j\in \mathbb N$ there exist $n_j \in \mathbb N$ and an operator $R_j\in \L(\ell_r^{n_j},\eell_s)$ such that $\|R_j\|_{\K_q}\leq$ $a_j$ and $\|R_j\|_{\K_p} > 2b_j$. Also, by Lemma~\ref{Lemma: On Codomain} and Lemma~\ref{lemma: mp size finite dimensional}, there exists $m_j \in \mathbb N$ such that the operator $\widetilde P_{m_j} \circ R_j\colon \ell_r^{n_j}\rightarrow \ell_s^{m_j}$ satisfies

$$
\|\widetilde P_{m_j} \circ R_j\|_{\K_q}\leq a_j \quad \mbox{and} \quad \|\widetilde P_{m_j} \circ R_j\|_{\K_p} > b_j.
$$ 

Since $E$ contains uniformly complemented copies of $\ell_r^n$, we can find $\lambda \ge 1$ so that there exist a subspace $E_{n_j} \subset E$, a projection $P_{n_j}\colon E \rightarrow E_{n_j}$ and an isomorphism $U_{n_j}\colon E_{n_j} \rightarrow \ell_{r}^{n_j}$ such that $B_{\ell_{r}^{n_j}} \subset U_{n_j}(B_{E_{n_j}}) \subset \lambda B_{\ell_{r}^{n_j}}$. On the other hand, if $\ell_s^n$ embed uniformly to $F$, then there exists an injective operator $T_{m_j}\colon \ell_s^{m_j}\rightarrow F$ with $\|T_{m_j}\|\|T^{-1}_{m_j}\|<\widetilde \lambda$ for some $\widetilde \lambda>1$. We may suppose that $\|T_{m_j}\|=1$.

Define a sequence of operators $V_j\colon E\rightarrow F$ as
$$
V_j=T_{m_j} \circ \widetilde P_{m_j} \circ R_{n_j} \circ U_{n_j} \circ P_{n_j}.
$$ 
Moreover, if $F$ contains uniformly complemented copies of $\ell_s^n$, we may also find linear operators $Q_j\colon F\rightarrow F$ with $\|Q_j\|\leq \widetilde \lambda$ such that $Q_j V_j=V_j$ and $Q_j V_{\widetilde j}=0$ if $j\neq \widetilde j$. This is because we have complemented copies of $\ell_s^n$ of any dimension $n\in \mathbb N$ as many times we need.
As $P_{n_j}$ is a finite rank operator, so is $V_j$. The estimation of the $\K_q$ norm of $V_j$ is easy:
\begin{equation}\label{Eq_1}
\|V_j\|_{\K_q}\leq \|T_{m_j}\| \|\widetilde P_{m_j} \circ R_{j}\|_{\K_q} \|U_{n_j}\| \|P_{n_j}\|<\lambda a_j.
\end{equation}

To compute the $\K_p$-norm of $V_j$, note that
$$
V_j(B_E)=T_{m_j} \circ \widetilde P_{m_j} \circ R_j \circ U_{n_j}(B_{E_{n_j}})\supset  T_{m_j} \circ \widetilde  P_{m_j} \circ R_j(B_{\ell_r^{n_j}})=T_{m_j}\left(\widetilde P_{m_j} \circ R_j(B_{\ell_r^{n_j}})\right).
$$
Then, since for any $p$-compact set $K$ in $\ell_s^{m_j}$ we have the inequalities
$$\m_p(K;\ell_s^{m_j})=\m_p(T_{m_j}^{-1}\circ T_{m_j}(K);\ell_s^{m_j})\leq \|T_{m_j}^{-1}\| \m_p(T_{m_j}(K);F)<\widetilde \lambda \m_p(T_{m_j}(K);F),$$
it follows that
\begin{equation}\label{Eq_2}
\begin{array}{rl}
\|V_j\|_{\K_p}=&\d \m_p(V_j(B_E);F)\\
\geq& \d  \m_p(T_{m_j}\left(\widetilde P_{m_j} \circ R_j(B_{\ell_r^{n_j}})\right);F)\\
\geq & \d  \frac{1}{\widetilde{\lambda}} \m_p(\widetilde P_{m_j} \circ R_j(B_{\ell_r^{n_j}});\ell_s^{m_j})\\
=& \d  \frac{1}{\widetilde \lambda} \|\widetilde P_{m_j} \circ R_j\|_{\K_p}\\
\geq &\d  
\frac{ \ b_j }{\widetilde \lambda}.
\end{array}
\end{equation}

It follows from \eqref{Eq_1} and \eqref{Eq_2} that
$$\|V_j\|_{\K_p} / \|V_j\|_{\K_q} \geq \frac{b_j}{\lambda \cdot \widetilde \lambda a_j},$$
which, since $M$ was arbitrary, proves \eqref{eq_mainthm}. The proof of a) is complete since we are assuming that $\K_q(E,F) \setminus \K_p(E,F) \not= \emptyset$. To see b), it is enough to show that $\K_q(E,F) \setminus \K_p(E,F) \not= \emptyset$. For this, note that the operator $V=\sum_{j=1}^{\infty} V_j$ is $q$-compact and $$\|V\|_{\K_q}\leq\sum_{j=1}^{\infty}\|V_j\|_{\K_p}=\lambda \|(a_j)_j\|_{\ell_1}.$$ Also, $V$ is not $p$-compact. Indeed, note that for each $j\in \mathbb N$ we have $V_j=Q_jV$. Then, if $V$ were $p$-compact, by \eqref{Eq_2} we would have $$\frac{b_j}{\widetilde \lambda}\leq \|V_j\|_{\K_p}=\|Q_j V\|_{\K_p}\leq \widetilde \lambda \|V\|_{\K_p},$$ which cannot happen since $\lim\limits_{j\to +\infty} b_j=\infty$.
\end{proof}

Combining the above theorem with Theorem~\ref{Thm:subspace lr} we obtain the following corollary that will be used.

\begin{corollary}\label{Coro: Application}
Let $E$ be a Banach space and let $1\leq p<q< \infty$ and $1\leq r,s\leq \infty$. Suppose that $\K_q(\eell_r,\eell_s)\setminus \K_p(\eell_r,\eell_s)$ is non-empty and $E$ contains uniformly complemented copies of $\ell_r^n$. Then $\K_q(E,\eell_s) \setminus \K_p(E,\eell_s)$ is $\eell_s$-spaceable.
\end{corollary}

The classical Dvoretzki's theorem states that any infinite-dimensional Banach space contains uniformly copies of $\ell_2^n$. On the other hand, it follows from Propositions \ref{Prop: q=p-compact} and \ref{Prop: example q=p-comp s>2} that $\K_q(\ell_1, \ell_2) \setminus \K_p(\ell_1,\ell_2) \not= \emptyset$ whenever $q > \max\{2, p\}$. These two facts provide one more consequence of Theorem \ref{Thm: Main_Linear}:

\begin{corollary} \label{Cor.linear.case}
	Let $E$ and $F$ be infinite-dimensional Banach spaces and let $1 \leq p < q < \infty$. Suppose that $E$ contains uniformly complemented copies of $\ell_1^n$ and $q > 2$. Then $\K_q(E, F) \setminus \K_p(E,F)$ is empty or spaceable.
\end{corollary}

Before ending, we would like to compare our result with that of Hern\'andez, Ruiz and S\'anchez. Following \cite[Definition~3.1]{HeRuSan} we say that a Banach space E is {\it $\sigma$-reproducible} if there exists a sequence $(E_n)_n$ of complemented subspaces, where $P_n\colon E\rightarrow E_n$ is a bounded projection, such that each $E_n$ is isomorphic to $E$, $P_i\circ P_j=0$ if $i \neq j$, and for all $k\in \mathbb N$ the projections $\widetilde P_k=\sum_{n=1}^{k}P_n\colon E\rightarrow \bigoplus^k_{n=1} E_n$ are uniformly bounded. In \cite[Theorem~3.5]{HeRuSan} is stated that if $E$ or $F$ are $\sigma$-reproducible, then $\K_q(E,F) \setminus \K_p(E,F)$ is empty or spaceable.  
First note that any Banach space which contains an infinite-dimensional complemented indecomposable Banach space is not $\sigma$-reproducible (recall that an infinite-dimensional Banach space $E$ is {\it indecomposable} if there do not exist infinite-dimensional closed subspaces $F$ and $G$ of $E$ with $E=F\oplus G$).
Now take $G$  an infinite-dimensional indecomposable Banach space and consider the 2-fold symmetric tensor product of $G$, $E=\widehat{\bigotimes}^{2,s}_{\pi_s}G$. It follows from \cite[Corollary 4]{BLA} that $E$ contains a complemented copy of $G$. Hence, neither $E$ or $G$ are not $\sigma$-reproducible and we cannot apply \cite[Theorem~3.5]{HeRuSan} to prove the spaceability of $\K_q(E,G)\setminus \K_p(E,G)$. But, since by \cite[Proposition~1.9]{BFV} $G$ contains uniformly complemented copies of $\ell_1^n$, we can apply Corollary~\ref{Cor.linear.case} to show the spaceability of $\K_q(E,G)\setminus \K_p(E,G)$. We will use this fact in the next section.

\section{The non-linear cases}\label{Sec: non linear case}

We now turn into the study of spaceability in the spaces of non-linear functions which are determinated by $p$-compact sets. We start with $k$-homogeneous polynomials. A function $P$ between Banach spaces $E$ and $F$ is a $k$-homogeneous polynomial if there exists a (unique) symmetric $k$-linear operator $\widehat P\colon E \times \cdots \times E\rightarrow F$ such that $P(x)=\widehat P(x,\ldots,x)$. The space of all $k$-homogeneous polynomials from $E$ into $F$ is denoted by $\P(^kE,F)$ and it is a Banach space when endowed with the norm $\|P\|=\sup_{x\in B_E} \|P(x)\|$. One of the main tool which is used to study different types of polynomials is the linearization via the $k$-fold symmetric tensor product, that we now describe briefly. Given Banach spaces $E$ and $F$ and a $k$-homogeneous polynomial $P\colon E\rightarrow F$, there exists a unique linear operator $T_P \in \L(\widehat \bigotimes^{k,s}_{\pi_s} E, F)$ such that 
$$
P(x)=T_P(x\otimes x \otimes\cdots \otimes x),
$$ 
where $\widehat \bigotimes^{k,s}_{\pi_s} E$ is the complete $k$-symmetric tensor product of $E$ endowed with the $k$-symmetric projective tensor norm. With this identification, combined with the fact that $\|P\|=\|T_P\|$, we obtain that $\P(^kE,F)$ is isometrically isomorphic to $\L(\widehat \bigotimes^{k,s}_{\pi_s} E,F)$. For more of this we refer the reader to the book of Dineen \cite{Din}.

Aron, Maestre and Rueda \cite{AMR} introduced $p$-compact $k$-homogeneous polynomials, extending the concept of compact $k$-homogeneous polynomial introduced by Aron and Schottenloher \cite{AS}. For $1\leq p <\infty$,  a $k$-homogeneous polynomial $P\in \P(^kE,F)$ is $p$-compact (compact) if $P(B_E)$ is a $p$-compact (compact) set. The space of all $p$-compact  $k$-homogeneous polynomials from $E$ to $F$ is denoted by $\P_{\K_p}(^kE,F)$ and it is a Banach space endowed with the norm $\|P\|_{\P_{\K_p}}=\m_p(P(B_E);F)$. 

In \cite[Theorem~3.1]{AR} is stating that there is an isometric isomorphism between $\P_{\K_p}(^kE,F)$ and $\K_p(\widehat \bigotimes^{k,s}_{\pi_s} E,F)$, where the isomorphism is given by $P\mapsto T_P$. This result, together with Theorem~\ref{Thm: Main_Linear} leads us to study the geometry of $\widehat \bigotimes^{k,s}_{\pi_s} E$ and $F$ in order to study the spaceability in the spaces of $p$-compact $k$-homogeneous polynomials from $E$ to $F$.  As we said before, from \cite[Proposition~1.9]{BFV} we have that for every infinite dimensional space $E$, $\widehat \bigotimes^{k,s}_{\pi_s} E$ contains uniformly complemented copies of $\ell_1^n$. We now put all this together and we obtain the following result.

\begin{proposition}\label{Prop: Polynomial}
Let $E$ and $F$ be infinite-dimensional Banach spaces and let $1\leq p<q<\infty$ and $1\leq s\leq \infty$. Suppose that $\K_q(\ell_1,\eell_s)\setminus \K_p(\ell_1,\eell_s)$ is non-empty and $\ell_s^n$ embed uniformly to $F$. Then $\P_{\K_q}(^kE,F)\setminus \P_{\K_p}(^kE,F)$ is empty or spaceable whenever $k \geq 2$.
\end{proposition}

\begin{proof}
The proof is basically posted in the paragraph above. Fix $E$ and $F$ Banach spaces. Since by \cite[Theorem~3.1]{AR}, $\P_{\K_p}(^kE,F)$ is isometrically isomorphic to $\K_p(\widehat \bigotimes^{k,s}_{\pi_s} E,F)$, we have that $\P_{\K_q}(^kE,F)\setminus \P_{\K_p}(^kE,F)$ is spaceable iwef and only if $\K_q(\widehat \bigotimes^{k,s}_{\pi_s} E,F)\setminus \K_p(\widehat \bigotimes^{k,s}_{\pi_s} E,F)$ is spaceable. Since by \cite[Proposition~1.9]{BFV}, $\widehat \bigotimes^{k,s}_{\pi_s} E$ contains uniformly complemented copies of $\ell_1^n$ and we are assuming that $\ell_s^n$ embed uniformly to $F$, and $\K_q(\ell_1,\eell_s)\setminus \K_p(\ell_1,\eell_s)$ is non-empty, an application of Theorem~\ref{Thm: Main_Linear} completes the proof.
\end{proof}

The study of spaceability in case of Lipschitz functions which are determined by $p$-compact sets follows the same pattern of what was done with $k$-homogeneous polynomials. For Banach spaces $E$ and $F$, denote by $\Lip_0(E,F)$ the Banach space of all Lipschitz functions from $E$ to $F$ which maps $0$ to $0$, endowed with the norm $\Lip(f)=\inf_{x\neq y} \frac{\|f(x)-f(y)\|}{\|x-y\|}$. When $F= \mathbb{K}$, $\Lip_0(E,\mathbb{K})=E^{\#}$ is the Lipschitz dual of $E$. Denote by $\delta\colon E\rightarrow \Lip_0(E^{\#})'$ the Dirac map given by $\delta(x) (f)=f(x)$. The Arens--Eells space of $E$, introduced by Arens and Eells \cite{AE} is defined as the norm-closed linear span of $\delta(E)$ in the space $(E^{\#})'$ and it is denoted by $\Ae(E)$. For more of this space we refer to the reader to the monograph of Godefroy \cite{God}, the book of Weaver \cite{Wea} and the reference therein. The Arens--Eells space can be used to linearize Lipschitz functions. Given a function $f\in \Lip_0(E,F)$ there exists a (unique) linear function $L_f\in \L(\Ae(E),F)$ such that $f(x)=L_f\circ \delta(x)$ (see, e.g. \cite[Theorem~3.6]{Wea}). Since also $\Lip(f)=\|L_f\|$, the map $f\mapsto L_f$ establishes an isometric isomorphism from $\Lip_0(E,F)$ onto $\L(\Ae(E),F)$.

For $1\leq p<\infty$, $p$-compact Lipschitz functions were introduced in \cite{AchDahTur} for Banach spaces $E$ and $F$ as those $f \in \Lip_0(E,F)$ such that the set $\{\frac{f(x)-f(y)}{\|x-y\|} \colon x\neq y\}\subset F$ is $p$-compact. The space of all Lipschitz $p$-compact functions from $E$ to $F$ is denoted by $\Lip_{0, \K_p}(E,F)$  and is a Banach space endowed with the norm $\|f\|_{\Lip_{0, \K_p}}=\m_p(\{\frac{f(x)-f(y)}{\|x-y\|} \colon x\neq y\}; F)$. From \cite[Proposition~3.2]{AchDahTur} we have $\Lip_{0, \K_p}(E,F)$ is isometrically isomorphic to $\K_p(\Ae(E),F)$ via the map $f\mapsto L_f$. So, again, to study the spaceability in the space of $p$-compact Lipschitz function from $E$ to $F$ it is enough to study the geometry of the space $\Ae(E)$ and $F$ and use Theorem~\ref{Thm: Main_Linear}. As in the case of $k$-homogeneous polynomials, we have the following result.

\begin{proposition}\label{Prop: Lipschitz}
Let $E$ and $F$ be infinite-dimensional Banach spaces, $1\leq p<q<\infty$ and $1\leq s\leq \infty$. Suppose that $\K_q(\ell_1,\eell_s)\setminus \K_p(\ell_1,\eell_s)$ is non-empty and $\ell_s^n$ embed uniformly to $F$. Then $\Lip_{0,\K_q}(E,F)\setminus \Lip_{0,\K_p}(E,F)$ is empty or spaceable.
\end{proposition}

\begin{proof}
Fix $E$ and $F$ Banach spaces. Since by \cite[Proposition~3.2]{AchDahTur} $\Lip_{0,\K_p}(E,F)$ is isometrically isomorphic to $\K_p(\Ae(E),F)$, we have that $\Lip_{0,\K_q}(E,F)\setminus \Lip_{0,\K_p}(E,F)$ is spaceable if and only if $\K_q(\Ae(E),F)\setminus \K_p(\Ae(E),F)$ is spaceable. Since by \cite[Theorem 1.1]{CDW}, there exists a complemented subspace of $\Ae(E)$ which is isomorphic to $\ell_1$ (hence $\Ae(E)$ contains uniformly copies of $\ell_1^n$) and we are assuming that $\ell_s^n$ embed uniformly to $F$ and $\K_q(\ell_1,\eell_s)\setminus \K_p(\ell_1,\eell_s)$ is non-empty, an application of Theorem~\ref{Thm: Main_Linear} completes the proof.
\end{proof}

As a first application of Proposition~\ref{Prop: Polynomial} and Proposition~\ref{Prop: Lipschitz}, a direct application of of Dvoretski theorem  gives the following result
\begin{proposition}\label{Prop: Poly and Lipschitz general}
Let $E$ and $F$ be infinite-dimensional Banach spaces and let $1 \leq p < q < \infty$. If $q > 2$, then
\begin{enumerate}[\upshape a)]
\item  $\P_{\K_q}(^kE,F) \setminus \P_{\K_p}(^kE,F)$  is empty or spaceable, whenever $k \geq 2$.
\item $\Lip_{0,\K_q}(E,F)\setminus \Lip_{0,\K_p}(E,F)$ is empty or spaceable.
\end{enumerate}
\end{proposition}

As we did in Corollary~\ref{Coro: Application}, when the codomain is $\eell_s$, thanks to Theorem~\ref{Thm:subspace lr} we can say more about the spaceability of the sets $\Lip_{0,\K_q}(E,\eell_s)\setminus \Lip_{0,\K_p}(E,\eell_s)$ and $\P_{\K_q}(^kE,\eell_s)\setminus \P_{\K_p}(^kE,\eell_s)$.

\begin{corollary}\label{Coro: Application to nonlinear}
Let $E$ be an infinite-dimensional Banach  space, $1\leq p<q<\infty$ and $1\leq s\leq \infty$. If there exists a $q$-compact set in $\eell_s$ which is not $p$-compact, then 

\begin{enumerate}[\upshape a)]
\item $\P_{\K_q}(^kE,\eell_s)\setminus \P_{\K_p}(^kE,\eell_s)$ is $\eell_s$-spaceable, whenever $k \geq 2$.
\item $\Lip_{0,\K_q}(E,\eell_s)\setminus \Lip_{0,\K_p}(E,\eell_s)$ is $\eell_s$-spaceable.
\end{enumerate}
\end{corollary}
\begin{proof}
By Proposition~\ref{Prop: q=p-compact}, if there exists a $q$-compact set in $\eell_s$ which is not $p$-compact, then $\K_q(\ell_1,\eell_s)\setminus \K_p(\ell_1,\eell_s)$ is non-empty. The proof follows by applying Corollary~\ref{Coro: Application} and procceding as in Proposition~\ref{Prop: Polynomial} and Proposition~\ref{Prop: Lipschitz} and Proposition~\ref{Prop: Lipschitz}.
\end{proof}

In Propositions~\ref{Prop: example q=p-comp s>2} and \ref{Prop: example q=p-comp s<2} we saw in which cases $q$-compact and $p$-compact sets coincide in $\eell_s$. From this, we can deduce in which cases $\P_{\K_q}(^kE,\eell_s)\setminus \P_{\K_p}(^kE,\eell_s)$ and $\Lip_{0,\K_q}(E,\eell_s)\setminus \Lip_{0,\K_p}(E,\eell_s)$ is spaceable or not. The proof is straightforward.

\begin{proposition}\label{Prop: spaceability in poly}
Let $1\leq p<q<\infty$, $E$ a Banach space and $1\leq s<\infty$. Then
\begin{enumerate}
\item If $q\leq \min\{2,s\}$, then $\P_{\K_q}(^kE,\eell_s)\setminus \P_{\K_p}(^kE,\eell_s)$ ($k\geq 2$) and $\Lip_{0,\K_q}(E,\eell_s)\setminus \Lip_{0,\K_p}(E,\eell_s)$ are empty.
\item If $q>\min\{2,s\}$, then $\P_{\K_q}(^kE,\eell_s)\setminus \P_{\K_p}(^kE,\eell_s)$ ($k\geq 2$) and $\Lip_{0,\K_q}(E,\eell_s)\setminus \Lip_{0,\K_p}(E,\eell_s)$ are $\eell_s$-spaceable.
\end{enumerate}
\end{proposition}

Also, using Proposition~\ref{Prop: r>2,s<2} and Proposition~\ref{Prop: r>2,s>2}  we may compare in which cases $\P_{\K_q}(^k\eell_r,\eell_s)\setminus \P_{\K_p}(^k\eell_r,\eell_s)$ and $\Lip_{0,\K_q}(\eell_r,\eell_s)\setminus \Lip_{0,\K_p}(\eell_r,\eell_s)$ is spaceable and $\K_q(\eell_r,\eell_s)\setminus \K_p(\eell_r,\eell_s)$ is not, for $r\geq 2$. All this is resumed in the following

\begin{proposition}\label{Prop: Comparation} Let $1\leq p<q<\infty$ and $r\geq 2$.
\begin{enumerate}
\item For $s\leq 2$ and $p\geq 2$, $\K_q(\eell_r,\eell_s)\setminus \K_p(\eell_r,\eell_s)$ is empty, meanwhile $\P_{\K_q}(^k E,\eell_s)\setminus \P_{\K_p}(^kE,\eell_s)$ ($k\geq 2$) and $\Lip_{0,\K_q}(E,\eell_s)\setminus \Lip_{0,\K_p}(E,\eell_s)$ are $\eell_s$-spaceable.
\item For $s\geq 2$ and $q>2$, $\K_q(\eell_r,\eell_s)\setminus \K_p(\eell_r,\eell_s)$ is empty, meanwhile $\P_{\K_q}(^kE,\eell_s)\setminus \P_{\K_p}(^k E,\eell_s)$ ($k\geq 2$) and $\Lip_{0,\K_q}(E,\eell_s)\setminus \Lip_{0,\K_p}(E,\eell_s)$ are $\eell_s$-spaceable.
\end{enumerate}
\end{proposition}

We now study the spaceability in the case of holomorphic function which are determined by $p$-compact set. From now on, all Banach spaces we will consider are complex. Recall that, for Banach spaces $E$ and $F$, a function $f\colon E\rightarrow F$ is holomorphic if for every $x_0\in E$ there exists $\ep>0$ and a sequence of $n$-homogeneous polynomials $(P_n(f)(x_0))_n$ such that $f(x)=\sum_{n=0}^\infty P_n(f)(x_0) (x-x_0)$ uniformly for every $x \in x_0+\ep B_E$. The series is called the Taylor series of $f$ at $x_0$. We denote by $\H(E,F)$ the space of all holomorphic functions from $E$ to $F$. For a general background of holomorphic function we refer to the books of Dineen~\cite{Din} and Mujica \cite{Mu}. 

Following Aron and Schottenloher \cite{AS}, a  holomorphic function $f \in \H(E,F)$ is {\it compact} if for every $x_0 \in E$ there exists $\ep>0$ such that $f(x_0+\ep B_E)$ is a relative compact set of $F$. The set of all compact holomorphic mappings from $E$ to $F$ is denoted by $\H_\K(E,F)$. The authors, characterize compact holomorphic function as those such that for every $x_0 \in E$ and $n \in \mathbb N$ the $n$-homogeneous polynomials in the Taylor series expansion of $f$ at $x_0$ are compact \cite[Proposition~3.4]{AS}. 

Aron, Maestre and Rueda \cite{AMR} introduced $p$-compact holomorphic functions in a {\it verbatim} way. A holomorphic function from $E$ to $F$ is said to be $p$-compact, $1\leq p <\infty$ if for every $x_0 \in E$ there exists $\ep>0$ such that $f(x_0+\ep B_E)$ is a relative $p$-compact set of $F$. However, the behavior of compact and $p$-compact mappings is different. Despite of if $f$ is a $p$-compact holomorphic function then every $n$-homogeneous polynomial in the Taylor series expansion of any $x_0$ in the domain is $p$-compact \cite[Proposition~3]{AMR}, in \cite[Example~3.7]{LaTur1} it was shown that there exists a holomorphic function $f\in \H(\ell_1,\ell_p)$ such that for every $x_0 \in \ell_1$ and every $n \in \mathbb N$ the $n$-homogeneous polynomials of the Taylor series expansion of $f$ at $x_0$ are $p$-compact but $f$ is not $p$-compact. Moreover, for every $x_0\in \ell_1$ and $\ep >0$, $f(x_0+\ep B_{\ell_1})$ fails to be a $p$-compact  set. This example answered a question posted in \cite{AMR}. 

We denote by $\H_{\K_p}(E,F)$ the set of all $p$-compact holomorphic functions from $E$ to $F$ and by $H_{\K_p}(E,F)$ the set of all holomorphic functions $f$ such that for every $x_0 \in E$, all the $n$-homogeneous polynomials in the Taylor series expansion of $f$ at $x_0$ are $p$-compact. Note that in virtue of \cite[Proposition~3.4]{AS} and \cite[Proposition~3]{AMR}, we have the inclusions $\H_{\K_p}\subset H_{\K_p} \subset \H_{\K}$ and, as we said before, \cite[Example~3.7]{LaTur1} tells us that $H_{\K_p}(\ell_1,\ell_p) \setminus \H_{\K_p}(E,F)$ is non-empty. We will finish the article by showing that $H_{\K_p}(\ell_1,\ell_p) \setminus \H_{\K_p}(E,F)$ is spaceable when we consider the Nachbin topology, also known as the {\it ported topology}. This topology was introduced by Nachbin and is denoted by $\tau_{\omega}$. An useful charaterization of this topology can be found in \cite[Propostion~3.47]{Din}: In $\H(E,F)$, the Nachbin topology is generated by the seminorms
$$
p(f)=\sum_{n=0}^{\infty} \sup_{x \in K+\alpha_n B_E} \| P_{n}f(0) (x)\|
$$
where $K$ ranges over all compact and balanced subset of $E$ and $(\alpha_n)_n\in c_0$.  

It is worth mentioning that $(\H_{\K}(E,F),\tau_{\omega})$ is a complete locally convex space and, as a consequence of \cite[Proposition~16.10]{Chae}, it is not metrizable when $E$ is an infinite-dimensional Banach space. In particular, in this context we cannot apply \cite[Theorem 7.4.2]{AronBerPelSeo} or {\cite[Theorem~3.3]{KitTim}.

We will consider
a partition of $\mathbb N$, $(L_j)_j$, such that every $L_j=\{n_1^j,n_2^j,\ldots\}$ has infinite elements and $L_j\cap L_k=\emptyset$ for $j\neq k$. Also, for $1\leq p <\infty$, we will use the operators $\iota_j\colon \ell_p\rightarrow \ell_p$ and $\pi_j\colon \ell_p\rightarrow \ell_p$ defined as
$$
\iota_j(x)=\sum_{i=1}^{\infty} e'_i(x) e_{n_i^j} \quad \mbox{and} \quad \pi_j(x)=\sum_{i=1}^{\infty} e'_{n_i^j}(x) e_i.
$$

\begin{lemma}\label{Lemma: last section} Let $1\leq p <\infty$ and take $K\subset \ell_p$ a $p$-compact set. For any sequence $(\beta_j)_j \in \ell_p$ the set 
$$
\widehat K=\{x \in \ell_p\colon x=\sum_{j=1}^{\infty} \beta_j \iota_j(y) \quad \mbox{with} \quad y \in K\}
$$
is $p$-compact and $\m_p(\widehat K;\ell_p)\leq \|(\beta_j)_j\|_{\ell_p} \m_p(K;\ell_p)$. 
\end{lemma}
\begin{proof}
As $K\subset \ell_p$ is $p$-compact, for $\ep>0$ there exists a sequence $(y_n)_n \in \ell_p(\ell_p)$ such that $K\subset \pco\{(y_n)_n\}$ and $\|(y_n)_n\|_{\ell_p(\ell_p)} \leq (1+\ep)\m_p(K;\ell_p)$. Fix $(\beta_n)_n\in \ell_p$ and, for each $n \in \mathbb N$, consider 
$$
z_n= \sum_{j=1}^{\infty}\beta_j \iota_j(y_n).
$$
Since we took a disjoint partition of $\mathbb N$, we have that $\|z_n\|_{\ell_p}= \|(\beta_j)_j\|_{\ell_p} \|y_n\|_{\ell_p}$. Thus, the sequence $(z_n)_n$ belongs to $\ell_p(\ell_p)$ and $\|(z_n)_n\|_{\ell_p(\ell_p)}=\|(\beta_j)_j\|_{\ell_p} \|(y_n)_n\|_{\ell_p(\ell_p)}$.

It only remains to see that $\widehat K \subset \pco\{(z_n)_{n}\}$. Take $x\in \widehat K$. Then, $x=\sum_{j=1}^{\infty} \beta_j \iota_j(y)$ for some $y \in K$. Thus, there exists a sequence $(\alpha_n)_n \in B_{\ell_{p'}}$ such that $y=\sum_{n=1}^{\infty} \alpha_n y_n$. Hence, 
$$
x=\sum_{j=1}^{\infty} \sum_{n=1}^{\infty} \beta_j\alpha_n \iota_j(y_n)=\sum_{n=1}^{\infty} \alpha_n \sum_{j=1}^{\infty} \beta_j \iota_j(y_n)=\sum_{n=1}^{\infty}\alpha_n z_n,
$$
and the proof follows.

\end{proof}
\begin{proposition}\label{Prop: Space_holo} For $1\leq p<\infty$, the set $H_{\K_p}(\ell_1,\ell_p) \setminus \H_{\K_p}(E,F)$ is spaceable as a subset of $(\H_{\K}(E,F),\tau_{\omega})$.
\end{proposition}
\begin{proof}
Take $f\in H_{\K_p}(\ell_1,\ell_p) \setminus \H_{\K_p}(E,F)$such that $P_1(f)(0) (e_1)=e'_1\otimes e_1$, which exist by \cite[Example~3.7]{LaTur1}, and define for each $j \in \mathbb N$ the functions $g_j=i_j \circ f$ and consider the vector space 
$$
W=span \{g_j\colon j \in \mathbb N\}.
$$ 
First we will see that $W\subset H_{\K_p}(\ell_1,\ell_p) \setminus \H_{\K_p}(\ell_1,\ell_p)\cup\{0\}$. Indeed, if $g\neq 0\in W$,  there is a sequence $(\beta_j)_j \in \mathbb C$ with finite nonzero elements such that $g=\sum_{j=1}^\infty \beta_j g_j$. Since for each $x \in \ell_1$ and each $n \in \mathbb N$, $P_ng(x)=\sum_{j=1}^{\infty} \iota_j \circ P_nf(x)$ is a finite sum, we conclude that $g_j \in H(\ell_1,\ell_p)$. If $g \in \H_{\K_p}(\ell_1,\ell_p)$, then take $j \in \mathbb N$ such that $\beta_j \neq 0$ and, since $f=\pi_{j} \circ g$, we obtain that $f\in \H_{\K_p}(\ell_1,\ell_p)$, which is a contradiction. Thus $W\subset H_{\K_p}(\ell_1,\ell_p) \setminus \H_{\K_p}(\ell_1,\ell_p)\cup\{0\}$.

To finish the proof, we will show that the $\tau_{\omega}$-closure of $W$ belongs to $H_{\K_p}(\ell_1,\ell_p) \setminus \H_{\K_p}(\ell_1,\ell_p)\cup\{0\}$. Fix $g\in \overline{W}^{\tau_{\omega}}$ and now we will obtain a description of $g$. There exists a net $(g_\alpha)_{\alpha \in \Gamma} \subset W$ such that $g_\alpha$ converges to $g$. For each $\alpha$, there exists a sequence $(\beta_j^{\alpha})_j \subset \mathbb C$ with finite nonzero elements such that $g_{\alpha}=\sum_{j=1}^{\infty} \beta_j^{\alpha} \iota_j\circ f$. As $g_{\alpha}$ converges to $g$, for each $j\in \mathbb N$, $\pi_j\circ g_{\alpha}=\beta_j^{\alpha} f$ $\tau_{\omega}$-converges to $\pi_j \circ g$.  In particular, for each $n \in \mathbb N$ and $x \in \ell_1$, the net $\beta_j^{\alpha} P_nf(x)$ converges in the usual polynomial norm to $\pi_j \circ P_n(g)(x)$. This implies that for each $j\in \mathbb N$ there exists $c_j \in \mathbb C$ such that $\beta_j^\alpha$ converges to $c_j$ as $\alpha\rightarrow +\infty$ and 
\begin{equation}\label{eq: pol g}
\pi_j \circ P_n(g)(x)=c_j P_nf(x)
\end{equation}
 for any $n\in \mathbb N$ and $x \in \ell_1$. 
In particular, since $P_1(f)(0) (e_1)=e_1$, we have that for any $j \in \mathbb N$, 
\begin{equation}\label{eq:c_j}
\pi_j \circ P_1(g)(0) (e_1)=c_j e_1.
\end{equation} On the other hand, as the partition of $\mathbb N$ is disjoint, we have that
\begin{equation}\label{eq: c_j2}
\|P_1(g)(0) (e_1)\|_{\ell_p}=\left(\sum_{j=1}^{\infty} \|\pi_j \circ P_1(g)(0) (e_1)\|_{\ell_p}^{p}\right)^p.
\end{equation}
Combining \eqref{eq:c_j} and \eqref{eq: c_j2} we obtain $(c_j)_j \in \ell_p$.

Since $(c_j)_j \in \ell_p$ we may define a function $\widehat g=\sum_{j=1}^{\infty} c_j \iota_j \circ f$. It is clear that $\widehat g\in \H(\ell_1,\ell_p)$, and for any $x \in \ell_1$, $n, j\in \mathbb N$, $\pi_j \circ P_n(\widehat g)(x)= c_j P_n(f) (x)$.  By \eqref{eq: pol g} and the uniqueness of the Taylor series expansion of a holomorphic function, we conclude that 
$$ 
g=\sum_{j=1}^{\infty} c_j \iota_j \circ f.
$$

If $g\neq 0$ is $p$-compact, then take $j_0 \in \mathbb N$ such that $c_{j_0}\neq 0$ and we obtain that $c_{j_0} f =\pi_{j_0} \circ g$ is $p$-compact, which cannot occur. Finally, to see that $g\in H_{\K_p}(\ell_1,\ell_p)$, fix $n\in \mathbb N$, and note that for $x \in \ell_1$
$$
P_ng(x)= \sum_{j=1}^{\infty} c_j \iota_j \circ P_nf(x).
$$
Since $P_nf(x)$ is a $p$-compact $n$-homogeneous polynomial, by Lemma~\ref{Lemma: last section} we have that $P_ng(x)(B_{\ell_1})$ is a $p$-compact set, and the proof follows.
\end{proof}

We finish by noting that By \cite[Proposition~1]{AMR}, every $p$-compact polynomial is, in particular, a $p$-compact holomorphic function. In virtue of this, the study the spaceability of $\H_{\K_q}\setminus \H_{\K_p}$ or $H_{\K_q}\setminus H_{\K_p}$ can be derive from Proposition~\ref{Prop: Polynomial}.

\subsection*{Acknowledgements} We would like to thank Dani Carando for introduce us the topic of this article.


\begin{thebibliography}{99}
	

\bibitem{AchDahTur} D. Achour, E. Dahia, P. Turco, \textit{Lipschitz p-compact mappings.} Monatsh. Math. 189 (2019), 595–609.


\bibitem{AE} R.F. Arens, J.Eells Jr., \textit{On embedding uniform and topological spaces.} Pacific J. Math. 6 (1956), 397--403. 

\bibitem{AronBerPelSeo} R. Aron, L. Bernal-Gonz\'{a}lez, D. Pellegrino, J. B. Seoane-Sep\'{u}lveda, \textsf{Lineability: The search for linearity in Mathematics.} Monographs and Research Notes in Mathematics, Chapman \& Hall/CRC, Boca Raton, FL, 2016.

\bibitem{ACGM_Trans}  R. Aron, E. \c{C}ali\c{s}kan, D. Garc\'{\i}a, M. Maestre, \textit{Behavior of holomorphic mappings on p-compact sets in a Banach space.} Trans. Amer. Math. Soc. 368 (2016), 4855--4871.

\bibitem{AGS} R. M. Aron, V. I. Gurariy, J. B. Seoane-Sep\'{u}lveda, \textit{Lineability and spaceability of sets of functions on $\mathbb{R}$.} Proc. Amer. Math. Soc. 133 (2005), 795--803.

\bibitem{AMR} R. Aron, M. Maestre, P. Rueda, \textit{$p$-compact holomorphic mappings.} Rev. R. Acad. Cienc. Exactas F\'is. Nat. Ser. A Math.  104 (2010), 353--364.

\bibitem{AR} R. Aron, P. Rueda, \textit{$p$-Compact homogeneous polynomials from an ideal point of view.} Contemp. Math., vol. 547, Amer. Math. Soc., (2011), 61--71.

\bibitem{AS} R. Aron, M. Schottenloher, \textit{Compact holomorphic mappings on Banach spaces and the approximation property.} J. Funct. Analysis 21 (1976), 7--30.


\bibitem{BerPelSeo} L. Bernal-Gonz\'{a}lez, D. Pellegrino, J. B. Seoane-Sep\'{u}lveda, \textit{Linear subsets of nonlinear sets in topological vector spaces.} Bull. Amer. Math. Soc. 51 (2014), 71--130.

\bibitem{BLA} F. Blasco, \textit{Complementation in spaces of symmetric tensor products and polynomials.} Studia Math., 123 (1997), 165--173.


\bibitem{BFV} F. Bombal, M. F\'ernandez-Unzueta, I. Villanueva, \textit{Local structure and copies of $c_0$ and $\ell_1$ in the tensor product of Banach spaces.} Bol. Soc. Mat. Mexicana 10 (2004), 195--202.
.

\bibitem{Chae} S.B. Chae, \textsf{Holomorphy and Calculus in Normed Spaces.} Monographs and Textbooks in Pure and Appl. Math.,
vol. 92, Marcel Dekker, 1985.

\bibitem{CDW} M. C\'uth, M. Doucha, P. Wojtaszczyk, \textit{On the structure of Lipschitz-free spaces.} Proc. Ame. Math. Soc. 144 (2016), 3833--3846.

\bibitem{DF} A. Defant, K. Floret, \textsf{Tensor norms and operators ideal.} North Holland Publishing Co., Amsterdam, 1993.

\bibitem{DelPinSer} J. M. Delgado, C. Pi\~{n}eiro, E. Serrano, \textit{Operators whose adjoints are quasi $p$-nuclear.} Studia Math. 197 (2010), 291--304.

\bibitem{Din} S. Dineen, \textsf{Complex Analysis on Infinite Dimensional Spaces.} Springer, 1999.


\bibitem{FouSwa} J. Fourie, J. Swart, \textit{Banach ideals of p-compact operators.} Manuscripta Math. 26 (1979), 349--362.

\bibitem{GaLaTur} D. Galicer, S. Lassalle, P. Turco, \textit{The ideal of $p$-compact operators: a tensor product approach.} Studia Math. 211 (2012), 269--286.

\bibitem {God} G. Godefroy, \textit{A survey on Lipschitz-free Banach spaces.} Commentationes Mathematicae 55  (2015), 89--118.

\bibitem{Gu91} V. I. Gurariy, \textit{Linear spaces composed of non-differentiable functions.} C. R. Acad. Bulg. Sci. 44 (1991), 13--16.

\bibitem{GurQua} V. I. Gurariy, L. Quarta, \textit{On lineability of sets of continuous functions.} J. Math. Anal. Appl. 294 (2004), 62--72.

\bibitem{HeRuSan} F. Hern\'andez, C. Ruiz, V. S\'anchez. \textit{Spaceability and operator ideals.} J. Math. Anal. Appl. 431 (2105), 1035--1044.

\bibitem{KitTim} D. Kitson, M. Timoney, \textit{Operator range and spaceability.} J. Math. Anal. Appl. 378 (2011), 680--686.

\bibitem{LaTur1} S. Lassalle, P. Turco, \textit{On $p$-compact mappings and the $p$-approximation properties.} J. Math. Anal. Appl. 389 (2012), 1204--1221.

\bibitem{LaTur2} S. Lassalle, P. Turco, \textit{The Banach ideal of $\A$-compact operators and related approximation properties.} J. Funct. Anal. 265 (2013), 2452--2464.

\bibitem{Mu}  J. Mujica, \textsf{Complex Analysis in Banach Spaces.} Math. Studies, vol. 120, North-Holland, Amsterdam, 1986.

\bibitem{Per} A. Persson, \textit{On some properties of p-nuclear and p-integral operators.} Studia Math. 33 (1969), 213--222.

\bibitem{Pietsch_p_summ} A. Pietsch, \textit{Absolutely $p$-summing operators in $\L_{r}$-spaces.} Bull. Soc. Math. Fr., Suppl., M\'em. 31--32 (1972), 285--315.

\bibitem{Pietsch_Book} A. Pietsch, \textsf{Operators Ideals.} Deutsch. Verlag Wiss., Berlin, 1978; North-Holland Publishing Company, Amsterdam, New York, Oxford, 1980.

\bibitem{Pietsch_p_compact} A. Pietsch, \textit{The ideal of $p$-compact operators and its maximal hull.} Proc. Amer. Math. Soc. 142 (2014), 519--530.

\bibitem{DelPin} C. Pi\~neiro, J. M. Delgado, \textit{$p$-convergent sequences and Banach spaces in which $p$-compact sets are $q$-compact.} Proc. Amer. Math. Soc. 139 (2011), 957--967.


\bibitem{SiKa} D. P. Sinha, A. K. Karn, \textit{Compact operators whose adjoints factor trough subspaces of  $\ell_p$.} Studia Math. 150 (2002), 17--33.

\bibitem{Tur} P. Turco, \textit{$\A$-compact mappings.} Rev. R. Acad. Cienc. Exactas F\'is. Nat. Ser. A Math. 110 (2015), 863--880.

\bibitem{Wea} N. Weaver, \textsf{Lipschitz Algebras.}  Hackensack, NJ World Scientific Publishing, 2018. 

\end{thebibliography}
\end{document}